\documentclass[11pt,a4paper]{article}% For Computer Modern Fonts, or,
\usepackage{latexsym,amssymb}
\usepackage{graphics,times,mathptm}
\usepackage{epsfig}
\newenvironment{keywords}{ \noindent {\small\bf Key Words}:}{ }
\newenvironment{Acknowledgement}{ \noindent {\bf Acknowledgement}.}{ }

\def\beq{\begin{equation}}
\def\eeq{\end{equation}}
\def\bea{\begin{eqnarray}}
\def\eea{\end{eqnarray}}
\def\beas{\begin{eqnarray*}}
\def\eeas{\end{eqnarray*}}

\newtheorem{theorem}{Theorem}

\newtheorem{property}{Property}

\begin{document}
%\begin{article}

\title{Generation of symmetric exponential sums}

%% Author name
\newcommand{\nms}{\normalsize}
\author{{\small \bf Yaroslav D. Sergeyev}\\ \\ [-2pt]
         \nms D.E.I.S. -- Universit\`a della Calabria, 87030 Rende (CS)
         -- Italy\\[-4pt]
         \nms and University of Nizhni Novgorod,\\[-4pt]
         \nms Gagarin Av., 23, Nizhni Novgorod -- Russia\\[-4pt]
         \nms {\tt (yaro@si.deis.unical.it)}
}

\date{}

\maketitle

\begin{abstract}
In this paper,  a new method for generation of infinite series of
symmetric identities written for exponential sums in real numbers
is proposed. Such systems have numerous applications in theory  of
numbers, chaos theory, algorithmic complexity, dynamic systems,
etc. Properties of generated identities are studied. Relations of
the introduced method for generation of symmetric exponential sums
to the Morse-Hedlund sequence and to the theory of magic squares
are established.
 \end{abstract}


\section{Introduction}
Let us consider two sets of numbers, $\{x_i: 1 \le i \le k\}$ and
$\{y_i: 1 \le i \le k\}$ such that the following  system of
identities is satisfied.
  \beq
 \left\{ \begin{array}{lcl}
      x^1_1  + x^1_2 + x^1_3  + ... + x^1_{k-1}   + x^1_k & = & y^1_1  + y^1_2  + y^1_3  + ... + y^1_{k-1}   +
     y^1_k \\
      x^2_1  + x^2_2 + x^2_3  + ... + x^2_{k-1}   + x^2_k & = & y^2_1  + y^2_2  + y^2_3  + ... + y^2_{k-1}   +
     y^2_k \\
      x^3_1  + x^3_2 + x^3_3  + ... + x^3_{k-1}   + x^3_k & = & y^3_1  + y^3_2  + y^3_3  + ... + y^3_{k-1}   +
     y^3_k \\
& \ldots \ldots \ldots & \\
      x^n_1  + x^n_2 + x^n_3  + ... + x^n_{k-1}   + x^n_k & = & y^n_1  + y^n_2  + y^n_3  + ... + y^n_{k-1}   +
     y^n_k
     \end{array}
                            \right.   \label{4}
 \eeq

Systems of the type  (\ref{4})  have a significant importance in
mathematics.  Such systems were widely studied in the number
theory  (see \cite{[1],[2],[3],[4]}). For example, the
Tarry-Escott problem (see, for example, \cite{Adler2,Lehmer})
deals with finding disjoint sets $\{x_i: 1 \le i \le k\}$ and
$\{y_i: 1 \le i \le k\}$ of integers satisfying (\ref{4}). The
systems (\ref{4}) have also a straightforward connection to the
Hilbert-Kamke problem for integers (see \cite{[5], [6], [16]})
being an extension of the famous Waring problem (see \cite{[7]}).

Identities  (\ref{4}) considered as a description of the process
of exponential sums expansion have close relations to the chaos
theory (see \cite{[14]}) and its applications, for instance, in
economics (see \cite{[15]}). Integer solutions to (\ref{4}) in the
context of integrals of piece-wise constant functions, emergent
calculations,  and integer code series  can be viewed (see
\cite{[8], [9]}) in terms of structural complexity of systems and
algorithms. Particularly, studying structural complexity of
multiextremal functions gives new tools for solving global
optimization problems being a new actively raising research area
(see \cite{Horst_and_Pardalos_(1995),[8],
[9],Strongin_and_Sergeyev}).

Together with the system (\ref{4})  the following system
(\ref{41}) can be considered.
  \beq
 \left\{ \begin{array}{lcl}
      x^1_1  + x^1_2 + x^1_3  + ... + x^1_{ \theta_{1}-1}   + x^1_{ \theta_1} & = & y^1_1  + y^1_2  + y^1_3  + ... + y^1_{ \theta_{1}-1}   +
     y^1_{ \theta_1} \\
      x^2_1  + x^2_2 + x^2_3  + ... + x^2_{ \theta_{2}-1}   + x^2_{ \theta_2} & = & y^2_1  + y^2_2  + y^2_3  + ... + y^2_{\theta_{2}-1}   +
     y^2_{ \theta_2} \\
      x^3_1  + x^3_2 + x^3_3  + ... + x^3_{ \theta_{3}-1}   + x^3_{ \theta_3} & = & y^3_1  + y^3_2  + y^3_3  + ... + y^3_{\theta_{3}-1}   +
     y^3_{ \theta_3} \\
& \ldots \ldots \ldots & \\
      x^n_1  + x^n_2 + x^n_3  + ... + x^n_{k-1}   + x^n_k & = & y^n_1  + y^n_2  + y^n_3  + ... + y^n_{k-1}   +
     y^n_k
     \end{array}
                            \right.   \label{41}
 \eeq
Such systems have interesting connections (see \cite{[9]})  to
algorithmic complexity (see \cite{[10],[11],[12]}) and period
doubling  in physical systems (see \cite{[13],[9]}). In this
case, the numbers
$\theta_{1},\theta_{2},\ldots,\theta_{n-1},\theta_{n}=k$ from
(\ref{4}) become
 \beq
  \theta_{2}=2\theta_{1},\hspace{5mm}
\theta_{3}=2\theta_{2}, \hspace{5mm} \ldots,
\hspace{5mm}\theta_{n}=2\theta_{n-1} \label{theta}
 \eeq
 and, for example, for
$\theta_{1}=2$ the following piramide is obtained.
  \beq
  \begin{array}{rcl}
      x^1_1  + x^1_2 & = & y^1_1  + y^1_2   \\
      x^2_1  + x^2_2 + x^2_3  +  x^2_4 & = & y^2_1  + y^2_2  + y^2_3  +   y^2_4 \\
      x^3_1  + x^3_2 + x^3_3  + ... + x^3_7   + x^3_8 & = & y^3_1  + y^3_2  + y^3_3  + ... + y^3_7   +
     y^3_8 \\
& \ldots \ldots \ldots & \\
      x^n_1  + x^n_2 + x^n_3  + ... + x^n_{k-1}   + x^n_{2^n} & = & y^n_1  + y^n_2  + y^n_3  + ... + y^n_{k-1}   +
     y^n_{2^n}
     \end{array}   \label{42}
 \eeq

In this paper, for any integer power $n$, a method for generation
of  two sets $\{x_i: 1 \le i \le k\}$ and  $\{y_i: 1 \le i \le
k\}$ of real numbers satisfying both (\ref{4}) and (\ref{41}) for
certain numbers $\theta_{1}\le \theta_{2}\le \ldots,
\theta_{n-1}\le \theta_{n}=k$ is proposed. Particularly,
piramidale systems of the type (\ref{41}), (\ref{theta}) are
considered. Properties of the introduced sets $\{x_i: 1 \le i \le
k\}$ and $\{y_i: 1 \le i \le k\}$ are studied. It is shown that
the sets found for a given $n$ are a basis for the sets
satisfying (\ref{4}) and (\ref{41}) with the powers $n+1, n+2,
\ldots$. In order to obtain solutions to a system having a higher
power $n+1$ it is necessary to add certain numbers to the sets
$\{x_i: 1 \le i \le k\}$ and $\{y_i: 1 \le i \le k\}$ found for
the power $n$. In such a way, by increasing $n$, it is possible to
speak about sequences $\{x_i\}$ and $\{y_i\}$ satisfying a
sequence of systems (\ref{4}) and (\ref{41}).

The next section contains the main results. A few examples   are
given in the last section. It also presents relations of the
introduced method for generation symmetric exponential sums to the
Morse-Hedlund sequences (see \cite{Morse}) being a particular
case of the Prouhet sequences (see \cite{Prouhet}) and to the
theory of magic cubes developed in \cite{Adler1,Adler2}.

\section{Main results}

In this section, two sets of real numbers $\{x_i: 1 \le i \le k\}$
and $\{y_i: 1 \le i \le k\}$ satisfying (\ref{4}) and (\ref{41})
for any given power $n$ are constructed\footnote{An alternative
description of these sets proposed by the unknown referee and is
given in the Appendix. }. Let us start by choosing four real
numbers $a, b, c,$ and $d$ satisfying the following simple
equality
 \beq
     a + b = c + d.   \label{1}
\eeq

 Given numbers $a, b, c, d$ from (\ref{1}) and any real
numbers $k_1, k_2,\ldots, k_n,n \ge 2,$ let us introduce by
recursion functions $L^j(k_1, k_2,\ldots, k_{n-1})$ and $R^j(k_1,
k_2,\ldots, k_{n-1})$, where  $j=n-1, n$. The functions
$L^1(k_1), R^1(k_1), L^2(k_1), R^2(k_1)$ and $ L^2(k_1,k_2),
R^2(k_1,k_2)$ are defined as follows:
 \beq
     L^1(k_1) = (a+k_1 )^1  + (b+k_1 )^1  = x^1_1  + x^1_2 ,
\label{5}
 \eeq
where $x^1_1  = a+k_1$ and $ x^1_2  = b+k_1 ,$
 \beq
     R^1(k_1) = (c+k_1 )^1  + (d+k_1 )^1  = y^1_1  + y^1_2 ,
     \label{6}
 \eeq
where $y^1_1  = c+k_1 $ and $ y^1_2  = d+k_1 ,$
\[
     L^2(k_1) = (a+k_1 )^2  + (b+k_1 )^2  = x^2_1  + x^2_2 ,
\]
\[
     R^2(k_1) = (c+k_1 )^2  + (d+k_1 )^2  = y^2_1  + y^2_2.
\]
Then, setting $x_3= c+k_1 +k_2 , x_4= c+k_1 +k_2, $ and
     $R^2 (k_1 +k_2 ) = x^2_3  + x^2_4$
we define
\[
     L^2 (k_1,k_2) = L^2 (k_1) + R^2 (k_1 +k_2) = x^2_1  + x^2_2  + x^2_3  +
     x^2_4=
\]
\[
                (a+k_1 )^2  + (b+k_1 )^2  + (c+k_1 +k_2)^2  + (d+k_1
                +k_2)^2.
\]
 Analogously we introduce $y_3  = a+k_1 +k_2, y_4  = b+k_1 +k_2,$
and $L^2 (k_1 +k_2) = y^2_3  + y^2_4$ to define
\[
     R^2 (k_1 , k_2) = R^2 (k_1) + L^2 (k_1 +k_2) = y^2_1  + y^2_2  + y^2_3  + y^2_4 =
\]
\[
                (c+k_1)^2 + (d+k_1)^2  + (a+k_1 +k_2)^2  + (b+k_1 +k_2)^2.
\]

Now, using already  determined  functions,  we  can  introduce by
recursion functions  $L^j(k_1, k_2,\ldots, k_{n-1})$  where
$j=n-1, n$:
 \beq
  L^n(k_1, k_2,\ldots, k_{n-1},k_{n}) = L^n(k_1, k_2,\ldots, k_{n-1})
  +R^n(k_1+ k_n, k_2,\ldots, k_{n-1}) = \sum^{2^n}_{i=1} x^n_i,
     \label{7}
 \eeq
 where
\beq
   L^n(k_1, k_2,\ldots, k_{n-1}) = \sum^{2^{n-1}}_{i=1} x^n_i,
     \label{8}
 \eeq
and the numbers $x_i , 1 \le i \le 2^{n-1} ,$ have already been
determined previously for powers  $1,\ldots,n-1$ and the
remaining numbers $x_i , 2^{n-1} < i \le 2^{n} ,$ are calculated
as follows
 \beq
     x_i  = y_j  + k_n , \hspace{5mm} i = 2^{n-1}   +j, \hspace{3mm}1 \le i \le 2^{n-1}  .
     \label{9}
 \eeq
The functions $R^j(k_1, k_2,\ldots, k_{n-1}),j=n-1, n,$ are
defined by a complete analogy:
 \beq
  R^n(k_1, k_2,\ldots, k_{n-1},k_{n}) = R^n(k_1, k_2,\ldots, k_{n-1})
  +L^n(k_1+ k_n, k_2,\ldots, k_{n-1}) = \sum^{2^n}_{i=1} y^n_i,
      \label{10}
 \eeq
\beq
   R^n(k_1, k_2,\ldots, k_{n-1}) = \sum^{2^{n-1}}_{i=1} y^n_i,
     \label{11}
 \eeq
 where $y_i , 1 \le i \le 2^{n-1} ,$ have already been determined  and
      \beq
     y_i  = x_j  + k_n , \hspace{5mm} i = 2^{n-1}   +j, \hspace{3mm}1 \le i \le 2^{n-1}  .
     \label{12}
 \eeq
Let us write for illustration the functions $L^3(k_1,k_2 ,k_3 ) $
and $ R^3(k_1,k_2 ,k_3 )$ in their explicit form.
\[
     L^3(k_1,k_2 ,k_3 ) = L^3(k_1,k_2  ) + R^3(k_1+k_2 ,k_3 ) =
\]
\[
 = (a+k_1 )^3     + (b+k_1 )^3     + (c+k_1+k_2 )^3     + (d+k_1+k_2 )^3  +
\]
\[
 + (c+k_1+k_3 )^3  + (d+k_1+k_3)^3  + (a+k_1+k_2 +k_3 )^3  + (b+k_1+k_2 +k_3)^3,
\]
\[
     R^3(k_1,k_2 ,k_3 ) = R^3(k_1,k_2  ) + L^3(k_1+k_2 ,k_3 ) =
 \]
\[
    = (c+k_1 )^3     + (d+k_1 )^3     + (a+k_1+k_2 )^3     + (b+k_1+k_2)^3  +
 \]
\[
    + (a+k_1+k_3 )^3  + (b+k_1+k_3 )^3  + (c+k_1+k_2 +k_3)^3  + (d+k_1+k_2 +k_3)^3 .
    \]
Now we are ready to formulate the first theorem.
 \begin{theorem}
 If $L^1 (0) = R^1 (0)$, i.e., (\ref{1}) holds, then for all $n \ge
 1$ and for
all real numbers $k_1, k_2,\ldots, k_{n}$  it follows
 \beq
      L^n(k_1, k_2,\ldots, k_{n}) =  R^n(k_1, k_2,\ldots, k_{n}).
      \label{13}
 \eeq
     \label{t1}
\end{theorem}
   \textbf{Proof.}  Theorem is proved  by induction. For $n =  1$  we
   obtain from
(\ref{1}), (\ref{5}), and (\ref{6}) that $L^1 (k_1 ) = R^1 (k_1
)$.

Suppose now that (\ref{13}) holds for $i$ such that $1 \le i \le
n-1$ and prove that this assumption implies truth of (\ref{13})
for $i = n$.  To proceed  we  need some designations. We denote
by $card \{ \sum_{m}k_m  \}$ the number of items $k_m,1 < m \le
n,$ in the sum $ \sum_{m}k_m$. We also introduce sets $S_{ij}  ,1
\le i < n, 0 \le j < i,$ as follows
\[
     S_{i0}   = \{ 0 \},
\]
\beq
     S_{ij}   = \{ s_{ij}  : s_{ij}   = \sum_{m}k_m ,
      card \{ \sum_{m}k_m  \} = j, n-i<m<n \},\hspace{3mm} 0<j<i.
      \label{14}
 \eeq
 Thus, a set $S_{ij}$   contains as elements all possible sums $\sum_{m}k_m ,
n-i<m<n,$ having $j$ items.

Now we are ready to calculate the left part of (\ref{13}), i.e.,
$L^n(k_1, k_2,\ldots, k_{n})$. By using (\ref{7}) and (\ref{9})
we obtain
 \beq
     L^n(k_1, k_2,\ldots,k_{n-1}, k_{n}) = L^n(k_1, k_2,\ldots, k_{n-1})
      + \sum^{2^{n-1}}_{j=1} (y_j+k_{n})^n.
           \label{15}
 \eeq
By using Newton's binomial formula for every item
$(y_j+k_{n})^n,  1 \le j \le 2^{n-1},$  in (\ref{15}) and then
applying (\ref{11}) and designations (\ref{14}) we have
\[
     L^n(k_1, k_2,\ldots,k_{n-1}, k_{n})  = L^n(k_1, k_2,\ldots, k_{n-1})
      + R^n(k_1, k_2,\ldots, k_{n-1}) +
\]
 \beq
   \sum^{n-1}_{i=1} k^i_n   \left(  \begin{array}{c} n \\ i \end{array} \right)
  [ \sum^{i-1}_{j=0} \hspace{2mm} \sum_{s_{ij} \in S_{ij}} R^{n-i}(k_1+s_{ij}, k_2,k_3,\ldots, k_{n-i}) ]
   + 2^{n-1}k^{n}_{n} ,
             \label{16}
 \eeq
 where $ \left(  \begin{array}{c} n \\ i \end{array} \right)$ are binomial coefficients.

We calculate $R^n(k_1, k_2,\ldots,k_{n-1}, k_{n})$ by a complete
analogy. It follows from (\ref{10}) and (\ref{12}) that
 \beq
       R^n(k_1, k_2,\ldots,k_{n-1}, k_{n}) = R^n(k_1, k_2,\ldots, k_{n-1})
      + \sum^{2^{n-1}}_{j=1} (x_j+k_{n})^n.
                   \label{17}
 \eeq
Applying again in (\ref{17}) Newton's  binomial   formula   for
items $(x_j+k_{n})^n,  1 \le j \le 2^{n-1},$  and then using
(\ref{8}) and (\ref{14}) we can write
\[
    R^n(k_1, k_2,\ldots,k_{n-1}, k_{n})  = R^n(k_1, k_2,\ldots, k_{n-1})
      + L^n(k_1, k_2,\ldots, k_{n-1}) +
\]
 \beq
   \sum^{n-1}_{i=1} k^i_n   \left(  \begin{array}{c} n \\ i \end{array} \right)
  [ \sum^{i-1}_{j=0} \hspace{2mm} \sum_{s_{ij} \in S_{ij}} L^{n-i}(k_1+s_{ij}, k_2,k_3,\ldots, k_{n-i}) ]
   + 2^{n-1}k^{n}_{n} ,
             \label{18}
 \eeq
Since numbers $n-i$ are such that $1 \le n-i < n$ then, due to our
assumption, the following identities
\[
     L^{n-i}(k_1+s_{ij}, k_2,k_3,\ldots, k_{n-i}) = R^{n-i}(k_1+s_{ij}, k_2,k_3,\ldots, k_{n-i})
\]
hold for all $0 \leq j \leq i-1, 1 \leq i \leq n-1$. This result
considered together with (\ref{16}) and  (\ref{18}) completes
induction and proves the theorem. \rule{5pt}{5pt}

Now we are ready to obtain the main result  concerning systems of
the type (\ref{4}) for the case $k = 2^n$.

 \begin{theorem}
 If $L^1 (0) = R^1 (0)$, i.e., (\ref{1}) holds, then for all $n \ge
 1$ and for
all real numbers $k_1, k_2,\ldots, k_{n}$ the numbers $x_j, y_j,
1 \le j \le 2^{n},$  defined by (\ref{5}), (\ref{6}), (\ref{9}),
and (\ref{12}) satisfy (\ref{42}) and   the following system
 \beq
 \left\{ \begin{array}{lcl}
      x^1_1  + x^1_2 + x^1_3  + ... + x^1_{2^n-1}   + x^1_{2^n} &  = & y^1_1  + y^1_2  + y^1_3  + ... + y^1_{2^n-1}   +
     y^1_{2^n} \\
      x^2_1  + x^2_2 + x^2_3  + ... + x^2_{2^n-1}   + x^2_{2^n} &  = & y^2_1  + y^2_2  + y^2_3  + ... + y^2_{2^n-1}   +
     y^2_{2^n} \\
      x^3_1  + x^3_2 + x^3_3  + ... + x^3_{2^n-1}   + x^3_{2^n} &  = & y^3_1  + y^3_2  + y^3_3  + ... + y^3_{2^n-1}   +
     y^3_{2^n} \\
& {\scriptstyle \ldots \ldots } & \\
      x^n_1  + x^n_2 + x^n_3  + ... + x^n_{2^n-1}   + x^n_{2^n} &  = & y^n_1  + y^n_2  + y^n_3  + ... + y^n_{2^n-1}   +
     y^n_{2^n}
     \end{array}
                            \right.   \label{19}
 \eeq
     \label{t2}
\end{theorem}
\textbf{Proof.}    Let  us  introduce  auxiliary  functions
$M^i(k_1, k_2,\ldots,k_{n-1}, k_{n}), 1\leq i \leq n,$ as follows
\beq
     M^i(k_1, k_2,\ldots,k_{n-1}, k_{n}) = x^i_1  + x^i_2 + x^i_3  + ... + x^i_{2^n-1}   + x^i_{2^n}, 1\leq i \leq n.
\label{20}
 \eeq
From (\ref{14}) and (\ref{16}) we can write
\[
     M^i(k_1, k_2,\ldots,k_{n-1}, k_{n}) =
     \sum^{n-i}_{j=0} \hspace{2mm} \sum_{s_{ij} \in S_{ij}} L^{i}(k_1+s_{ij}, k_2,k_3,\ldots, k_{i}).
     \]
Due to  (\ref{13}) and (\ref{10}), we can rewrite this equality as
\[
    M^i(k_1, k_2,\ldots,k_{n-1}, k_{n}) =
        \sum^{n-i}_{j=0} \hspace{2mm} \sum_{s_{ij} \in S_{ij}}
         R^{i}(k_1+s_{ij}, k_2,k_3,\ldots, k_{i})=  \sum^{2^n}_{j=1} y^i_j.
     \]
where $1 \leq i \leq n$. Thus, the theorem has been completely
proved. \rule{5pt}{5pt}

Let us now establish some properties of the introduced numbers
$x_m, y_m, 1 \le m \le 2^{n},$ from (\ref{19}).

\begin{property}If there exists a number $i > 1$ such that $k_i   =  0$
then it follows for all $n \geq i$ that for every number $x_m , 1
\leq m \leq 2^n,$ there exists such a number $y_p$ that $x_m  =
y_p$.
     \label{p1}
 \end{property}
\textbf{Proof.} Truth of the fact follows immediately from
definitions (\ref{9}) and (\ref{12}). \rule{5pt}{5pt}

The next property establishes presence of symmetry not only with
respect to the sign `=' but also for inner subsumes of system
(\ref{4}).

\begin{property} For subsumes of the numbers $x_m, y_m,1
\leq m \leq 2^n,$  from (\ref{19}) the following identities hold
\[
       \sum^{(p+1)2^i}_{j=p2^i+1} x^i_j  =    \sum^{(p+1)2^i}_{j=p2^i+1}, y^i_j ,
         \hspace{5mm}    0 \leq p \leq 2^{n-i}   , \hspace{2mm} 1 \leq i \leq n.
\]
     \label{p2}
 \end{property}
\textbf{Proof.} This result is a straightforward consequence of
theorems~\ref{t1} and \ref{t2}. \rule{5pt}{5pt}

\begin{property}If the number $a + b + 2k_1  +  \sum^{n}_{i=2} k_i $ is an integer then
\[
    M^1(k_1, k_2,\ldots,k_{n-1}, k_{n})_{\mbox{mod} 2}    = 0.
\]
     \label{p3}
 \end{property}
\textbf{Proof.} This result is obtained from (\ref{5}) --
(\ref{12}), since the number $M^1(k_1, k_2,\ldots,k_{n-1},
k_{n})$ can be written in the form
\[
     M^1(k_1, k_2,\ldots,k_{n-1}, k_{n}) =   \sum^{2^{n}}_{i=1} x_i  =
      2^{n-1}    (a + b + 2k_1  + \sum^{n}_{i=2} k_i ). \hspace{2mm} \rule{5pt}{5pt}
\]

\begin{property}For  a fast  calculation of  the  numbers  $M^i(k_1, k_2,\ldots,k_{n-1},
k_{n}), 1 \leq i \leq n,$ from (\ref{20}) the following recursive
formula can be used
\[
     M^i(k_1, k_2,\ldots,k_{n-1},
k_{n}) = 2 M^i(k_1, k_2,\ldots,k_{n-1}) +
\]
 \beq
       \sum^{n-1}_{i=1} k^i_n   \left(  \begin{array}{c} n \\ i \end{array} \right)  M^{n-i}(k_1, k_2,\ldots,k_{n-1})
      +  2^{n-1}k^{n}_{n}.
  \label{22}
 \eeq
     \label{p4}
 \end{property}
\textbf{Proof.} Property is a straightforward corollary of
(\ref{18}) and (\ref{20}). \rule{5pt}{5pt}

Theorem 2 has been proved for numbers $k$ from (\ref{4}) equal to
$2^n$. However, it  is possible to generate solutions to (\ref{4})
for $k \neq 2^n$ by a suitable choice of the numbers $a,b,c,d,k_1,
k_2,\ldots, k_{n}$ in such a way that some numbers $x_j= y_i=0 $
and so can be excluded from consideration in (\ref{4}). It is also
possible to generate `asymmetric' sums in sense that the number of
items not equal to zero laying on the left from `=' does not
coincide with that number for the sum laying on the right (see
examples in the next section).

The following theorem extends     results established  above from
the case of two items in the basic identity (\ref{1}) to the case
with  $N$ items.

\begin{theorem} If, given parameters $k_1, k_2,\ldots, k_{n}$ in the recursive
process (\ref{5}) - (\ref{12})  the basic identity (\ref{1}) is
substituted by  the  identity
 \beq
      \sum^{N}_{j=1} a _j  =   \sum^{N}_{j=1} c _j
        \label{231}
 \eeq
then, the numbers $x_i, y_i$ generated by the scheme will be
solutions to systems (\ref{4}), (\ref{41}) for $k = 2^{n-1}N$.
     \label{t3}
\end{theorem}
\textbf{Proof.}  The theorem is proved by a complete analogy with
the basic case $N = 2$.  \rule{5pt}{5pt}

Properties similar to the properties established for the case
(\ref{1}) can be proved for (\ref{231}) too.

\section{Examples and relations to the Morse-Hedlund sequence}

%{\small
%\[
%(a+k_1)^1+(b+k_1)^1=(c+k_1)^1+(d+k_1)^1
%\]
% \beq
% (a+k_1)^2+(b+k_1)^2+(c+k_1+k_2)^2+(d+k_1+k_2)^2=(c+k_1)^2+(d+k_1)^2+(a+k_1+k_2)^2+(b+k_1+k_2)^2
% \label{2}
% \eeq
%\[
% (a+k_1)^3+(b+k_1)^3+(c+k_1+k_2)^3+(d+k_1+k_2)^3
% +(c+k_1+k_3)^3+(d+k_1+k_3)^3+(a+k_1+k_2+k_3)^3+(b+k_1+k_2+k_3)^3
% =(c+k_1)^3+(d+k_1)^3+(a+k_1+k_2)^3+(b+k_1+k_2)^3
% +(a+k_1+k_3)^3+(b+k_1+k_3)^3+(c+k_1+k_2+k_3)^3+(d+k_1+k_2+k_3)^3
%\]}
%\begin{figure}[h]
%  \begin{center}
%    \caption{An example of the
%problem (\ref{problem})--(\ref{const}).  \label{f1}}
% \hspace{-2cm}
%     \epsfig{ figure = Numeri_1.eps, width = 7in, height = 0.9in,  silent = yes }
%  \end{center}
%\end{figure}

Let  us  consider the first example for the case $n = 3$ with the
numbers
\[
     a = -1,\hspace{3mm} b = 2.1,\hspace{3mm} c = 3.4,\hspace{3mm} d = -2.3,
     \]
 \beq
      k_1  = 0,\hspace{3mm} k_2  = 1,\hspace{3mm} k_3  = -0.5.
                   \label{3}
 \eeq
Application of the introduced method to  the  data  (\ref{3})
gives the following piramide  (\ref{42}) of identities.
\[
\hspace{-1cm}  \begin{array}{rcl}
 {\scriptstyle (-1)^1+2.1^1} & {\scriptstyle =} &  {\scriptstyle 3.4^1+(-2.3)^1,}\\
 {\scriptstyle (-1)^2+2.1^2+4.4^2+(-1.3)^2}
  & {\scriptstyle =} &  {\scriptstyle 3.4^2+(-2.3)^2+0^2+3.1^2,} \\
 {\scriptstyle (-1)^3+2.1^3+4.4^3+(-1.3)^3 +2.9^3+(-2.8)^3+(-0.5)^3+2.6^3}
  & {\scriptstyle =} &  {\scriptstyle 3.4^3+(-2.3)^3+0^3+3.1^3+(-1.5)^3+1.6^3+3.9^3+(-1.8)^3.}
   \end{array}
\]
 The sums for the first, second and third identities are equal to
$1.1,\; 26.46,$ and $111.136,$ respectively.

Let us now use the data from   (\ref{3})  and construct the
corresponding systems of the type (\ref{4}). This operation leads
to two groups of identities presented below. The first group
consists of identities for $n=2$
\[
\hspace{-1cm}  \begin{array}{c}
 {\scriptstyle (-1)^1+2.1^1+4.4^1+(-1.3)^1}=
    {\scriptstyle 3.4^1+(-2.3)^1+0^1+3.1^1,}\\
 {\scriptstyle (-1)^2+2.1^2+4.4^2+(-1.3)^2}
  =  {\scriptstyle 3.4^2+(-2.3)^2+0^2+3.1^2}
   \end{array}
\]
where  sums in the first line are equal to $4.2$ and in the second
to $26.46$. The second group of identities is written for $n=3$.
\[
\hspace{-1cm}  \begin{array}{c}
 {\scriptstyle (-1)^1+2.1^1+4.4^1+(-1.3)^1
 +2.9^1+(-2.8)^1+(-0.5)^1+2.6^1}=
    {\scriptstyle 3.4^1+(-2.3)^1+0^1+3.1^1+(-1.5)^1+1.6^1+3.9^1+(-1.8)^1,}\\
 {\scriptstyle (-1)^2+2.1^2+4.4^2+(-1.3)^2 +2.9^2+(-2.8)^2+(-0.5)^2+2.6^2}
  =  {\scriptstyle 3.4^2+(-2.3)^2+0^2+3.1^2+(-1.5)^2+1.6^2+3.9^2+(-1.8)^2,} \\
 {\scriptstyle (-1)^3+2.1^3+4.4^3+(-1.3)^3 +2.9^3+(-2.8)^3+(-0.5)^3+2.6^3}
  =  {\scriptstyle 3.4^3+(-2.3)^3+0^3+3.1^3+(-1.5)^3+1.6^3+3.9^3+(-1.8)^3.}
   \end{array}
\]
 The sums for the
first,  second,  and third identities are equal to 6.4, 49.72,
and 111.136, respectively.

The second example is in integer numbers for $n = 4$. In order to
obtain repeated numbers $x_i, y_i$  and to demonstrate the nature
of a `hidden' symmetry appearance the following parameters
 \beq
     a = 1,\hspace{2mm} b = 3,\hspace{2mm} c = d = 2,\hspace{8mm}
      k_1  = 0,\hspace{2mm} k_2  = 1, \hspace{2mm}k_3  = -2, \hspace{2mm}k_4  = 3
       \label{23}
 \eeq
have been chosen. The pyramid from (\ref{42}) for these data is
shown below.
\[
\hspace{-2cm}  \begin{array}{rcl}
 {\scriptstyle 1^1+3^1} & {\scriptstyle =} &  {\scriptstyle 2^1+2^1,}\\
 {\scriptstyle 1^2+3^2+3^2+3^2}
  & {\scriptstyle =} &  {\scriptstyle 2^2+2^2+2^2+4^2,} \\
 {\scriptstyle 1^3+3^3+3^3+3^3 +0^3+0^3+0^3+2^3}
  & {\scriptstyle =} &  {\scriptstyle
  2^3+2^3+2^3+4^3+(-1)^3+1^3+1^3+1^3,}\\
 {\scriptstyle 1^4+3^4+3^4+3^4 +0^4+0^4+0^4+2^4 + 5^4+5^4+5^4+7^4 +2^4+4^4+4^4+4^4}
  & {\scriptstyle =} &  {\scriptstyle 2^4+2^4+2^4+4^4+(-1)^4+1^4+1^4+1^4
  +4^4+6^4+6^4+6^4 +3^4+3^4+3^4+5^4.}
   \end{array}
\]
The sums for the first, second, third, and forth identities are
equal, respectively, to 4, 28, 90, and 5320. In this record all
the numbers $x_i, y_i$ have been shown. Let us now eliminate from
the record items equal to  zero  and elements $x_i$ and $y_j$ such
that $x_i = y_j$. In this case the symmetry  is  present only in
an implicit form and is not seen from the record:
\[
\hspace{1cm}  \begin{array}{rcl}
                              1^1  + 3^1 & = & 2^1  + 2^1,\\
       1^2  + 3^2  + 3^2  + 3^2   & = &  2^2  + 2^2  + 2^2 + 4^2,\\
 3^3  + 3^3  + 3^3   & = &  2^3  + 2^3  + 4^3  + (-1)^3  + 1^3  + 1^3,\\
   5^4  + 5^4  + 7^4  + 4^4   & = &  2^4  + (-1)^4  + 1^4  + 1^4  + 6^4  + 6^4
   \end{array}
\]
Let us consider now how does system (\ref{19}) look like for
data  (\ref{23}).
\[
\hspace{-2cm}  \begin{array}{rcl}
 {\scriptstyle 1^1+3^1+3^1+3^1 +0^1+0^1+0^1+2^1 + 5^1+5^1+5^1+7^1 +2^1+4^1+4^1+4^1}
  & {\scriptstyle =} &  {\scriptstyle 2^1+2^1+2^1+4^1+(-1)^1+1^1+1^1+1^1
  +4^1+6^1+6^1+6^1 +3^1+3^1+3^1+5^1,}\\
 {\scriptstyle 1^2+3^2+3^2+3^2 +0^2+0^2+0^2+2^2 + 5^2+5^2+5^2+7^2 +2^2+4^2+4^2+4^2}
  & {\scriptstyle =} &  {\scriptstyle 2^2+2^2+2^2+4^2+(-1)^2+1^2+1^2+1^2
  +4^2+6^2+6^2+6^2 +3^2+3^2+3^2+5^2,} \\
 {\scriptstyle 1^3+3^3+3^3+3^3 +0^3+0^3+0^3+2^3 + 5^3+5^3+5^3+7^3 +2^3+4^3+4^3+4^3}
  & {\scriptstyle =} &  {\scriptstyle 2^3+2^3+2^3+4^3+(-1)^3+1^3+1^3+1^3
  +4^3+6^3+6^3+6^3 +3^3+3^3+3^3+5^3,}\\
 {\scriptstyle 1^4+3^4+3^4+3^4 +0^4+0^4+0^4+2^4 + 5^4+5^4+5^4+7^4 +2^4+4^4+4^4+4^4}
  & {\scriptstyle =} &  {\scriptstyle 2^4+2^4+2^4+4^4+(-1)^4+1^4+1^4+1^4
  +4^4+6^4+6^4+6^4 +3^4+3^4+3^4+5^4.}
   \end{array}
\]
 The sums are equal to 48, 208, 1008, and 5320.
Deleting the  explicit symmetry we obtain the following
identities where the symmetry is hidden.
\[
\hspace{1cm}  \begin{array}{rcl}
5^1  + 5^1  + 7^1  + 4^1   & = &  2^1  + (-1)^1  + 1^1  + 1^1  + 6^1  + 6^1,\\
5^2  + 5^2  + 7^2  + 4^2   & = &  2^2  + (-1)^2  + 1^2  + 1^2  + 6^2  + 6^2,\\
5^3  + 5^3  + 7^3  + 4^3   & = &  2^3  + (-1)^3  + 1^3  + 1^3  + 6^3  + 6^3,\\
5^4  + 5^4  + 7^4  + 4^4   & = &  2^4  + (-1)^4  + 1^4  + 1^4  +
6^4  + 6^4
\end{array}
\]
The corresponding sums are equal to 21, 115, 657, and 3907.

Let us now illustrate Theorem~\ref{t3} by an example having $N =
3$ in (\ref{231}) and shown in Table~\ref{tab1} for the power
$n=4$. In this example the following parameters have been used:
 \[
     a_1  = 1, \hspace{3mm}a_2  = 3, \hspace{3mm}a_3  = 7,\hspace{2cm}
       c_1  = 2,\hspace{3mm} c_2  = 4, \hspace{3mm}c_3  = 5,
 \]
 \[
     k_1  = 0,\hspace{3mm}   k_2 = -1, \hspace{3mm}k_3 = 1.3,\hspace{3mm} k_4 = -2.5.
 \]
The last line of Table~\ref{tab1} presents the values
\[
      \sum^{24}_{i=1} x^j_i,\hspace{8mm}\sum^{24}_{i=1} y^j_i,\hspace{1cm} 1 \le j \le
      4,
\]
showing that the numbers $x_i, y_i, 1 \le j \le 4,$ are solutions
to system (\ref{4}) for $n = 4$ and $k = 24$.

     Table~\ref{tab1} illustrates also Property~\ref{p2}.
In fact, it is split in boxes that contain groups of elements of
the sequences $\{x^j_i\}$ and $\{y^j_i\}, 1 \le j \le 4,$ having
equal sums. For example, for the power $2$, the box containing
the numbers $x^2_7, x^2_8, ..., x^2_{11}, x^2_{12}$ and   $y^2_7,
y^8_2 , ..., y^2_{11}, y^2_{12}$ shows that
\[
x^2_7  + x^2_8 +  ... + x^2_{11}   + x^2_{12} =   y^2_7  + y^2_8
+   ... + y^2_{11}   + y^2_{12}.
\]
To conclude this example let us note, that $x_6 =  y_2   =  4$.
Elimination of this couple of numbers from consideration gives us
a solution to (\ref{4}) for $n = 4$ and $k=23$. We can repeat this
procedure for pairs
$$x_8 = y_{12}   = 5.3, \hspace{8mm} x_{14}   = y_{18}   =
1.5,\hspace{8mm} x_{24}=y_{20} = 2.8$$
 arriving to a solution to (\ref{4}) for
$n = 4$ and $k = 20$.
    \begin{table}[t]\centering
   \caption{Example for $n=4$ and the starting identity $1+3+7=2+4+5$}
    \label{tab1}
 \vspace{5mm}
    {\small
\begin{tabular}{|r|rr|rr|rr|rr|}
\hline
 &&&&&&&& \\
&$x^1$&$y^1$&$x^2$&$y^2$&$x^3$&$y^3$&$x^4$&$y^4$\\
%&&&&&&&& \\
\hline
1&1&2&1&4&1&8&1&16\\
2&3&4&9&16&27&64&81&256\\
3&7&5&49&25&343&125&2401&625\\ \cline{2-3}
4&1&0&1&0&1&0&1&0\\
5&3&2&9&4&27&8&81&16\\
6&4&6&16&36&64&216&256&1296\\  \cline{2-5}
7&3,3&2,3&10,89&5,29&35,937&12,167&118,5921&27,9841\\
8&5,3&4,3&28,09&18,49&148,877&79,507&789,0481&341,8801\\
9&6,3&8,3&39,69&68,89&250,047&571,787&1575,2961&4745,8321\\
                                            \cline{2-3}
10&1,3&2,3&1,69&5,29&2,197&12,167&2,8561&27,9841\\
11&3,3&4,3&10,89&18,49&35,937&79,507&118,5921&341,8801\\
12&7,3&5,3&53,29&28,09&389,017&148,877&2839,8241&789,0481\\
                                             \cline{2-7}
13&-0,5&-1,5&0,25&2,25&-0,125&-3,375&0,0625&5,0625\\
14&1,5&0,5&2,25&0,25&3,375&0,125&5,0625&0,0625\\
15&2,5&4,5&6,25&20,25&15,625&91,125&39,0625&410,0625\\
                                              \cline{2-3}
16&-2,5&-1,5&6,25&2,25&-15,625&-3,375&39,0625&5,0625\\
17&-0,5&0,5&0,25&0,25&-0,125&0,125&0,0625&0,0625\\
18&3,5&1,5&12,25&2,25&42,875&3,375&150,0625&5,0625\\
                                               \cline{2-5}
19&-0,2&0,8&0,04&0,64&-0,008&0,512&0,0016&0,4096\\
20&1,8&2,8&3,24&7,84&5,832&21,952&10,4976&61,4656\\
21&5,8&3,8&33,64&14,44&195,112&54,872&1131,6496&208,5136\\
                                               \cline{2-3}
22&-0,2&-1,2&0,04&1,44&-0,008&-1,728&0,0016&2,0736\\
23&1,8&0,8&3,24&0,64&5,832&0,512&10,4976&0,4096\\
24&2,8&4,8&7,84&23,04&21,952&110,592&61,4656&530,8416\\
\hline
Sum&61,6&61,6&305,08&305,08&1599,724&1599,724&9712,6972&9712,6972\\
\hline
   \end{tabular}
}
\end{table}

\begin{center}
     \begin{table}[t]\centering
   \caption{The
Morse-Hedlund  sequence}
   \label{tab2}
 \vspace{5mm}
    {\small   \begin{tabular}{rrrrrrrrrrrrrrrrr} \hline
1&0&0&1&0&1&1&0&0&1&1&0&1&0&0&1&...\\
1&2&3&4&5&6&7&8&9&10&11&12&13&14&15&16&...\\
    \hline
   \end{tabular}
   }
\end{table}
\end{center}

The introduced method for generation of symmetric exponential
sums has  interesting relations to the Morse-Hedlund sequences
\cite{Morse} (being a particular case of the Prouhet sequences
\cite{Prouhet}) and to the theory of magic cubes developed in
\cite{Adler1,Adler2}. Let us write the first items of the
Morse-Hedlund  sequence (see Table~\ref{tab2}) and numerate its
items by natural numbers.

 It has been noticed (see, for example,
\cite{[8],Prouhet}) that this sequence splits the natural numbers
staying in the second line of Table~\ref{tab2} in two groups
satisfying (\ref{4}). The first group contains the numbers
corresponding to `$1$' in the Morse-Hedlund sequence and the
second group contains the numbers corresponding to `$0$'. Thus,
for the first $16$ numbers we obtain
\[
\hspace{-2mm} {\small   \begin{array}{rcl}
1^1  + 4^1  + 6^1  + 7^1 + 10^1  + 11^1  + 13^1  + 16^1  & = &  2^1  + 3^1  + 5^1  + 8^1  + 9^1  + 12^1  + 14^1  + 15^1,\\
1^2  + 4^2  + 6^2  + 7^2 + 10^2  + 11^2+ 13^2  + 16^2 & = &  2^2  + 3^2  + 5^2  + 8^2  + 9^2  + 12^2+ 14^2  + 15^2,\\
1^3  + 4^3  + 6^3  + 7^3 + 10^3  + 11^3+ 13^3  + 16^3 & = &  2^3
+ 3^3  + 5^3  + 8^3 + 9^3+ 12^3+ 14^3  + 15^3.
\end{array}       }
\]
It can be seen immediately that these identities are a particular
case of our results with
 \beq
     a = 1,\hspace{3mm} b = 4,\hspace{3mm} c = 2,\hspace{3mm} d = 3, \hspace{1cm}
     k_1  = 0, \hspace{3mm}k_2  = 2^2,\hspace{3mm} k_3  = 2^3.    \label{27}
\eeq
 Naturally, by taking  $k_n  = 2^n$ it is possible to generate
Prouhet identities  for any $n$.

    \begin{table}
    \centering
   \caption{The half-completed magic square  based on the Morse-Hedlund sequence}
    \label{tab3}
 \vspace{5mm}
\begin{tabular}{|c|c|c|c|}
 \hline
 &  &  & \\
 \hspace{15mm}  &   2 & 3& \hspace{15mm} \\
  &&&\\    \hline
 &&&\\
 5& \hspace{15mm} & \hspace{15mm} &8 \\
  &&&\\  \hline
 &&&\\
 9&  &  &12\\
  &&&\\  \hline
 &&&\\
  &14&15&  \\
   &&&\\ \hline
   \end{tabular}
\end{table}
    \begin{table}
    \centering
   \caption{The completed magic square  based on the Morse-Hedlund sequence}
    \label{tab4}
 \vspace{5mm}
\begin{tabular}{|c|c|c|c|}
 \hline
 \hspace{15mm} & \hspace{15mm} & \hspace{15mm} & \hspace{15mm} \\
16  &   2 & 3& 13 \\
  &&&\\    \hline
 &&&\\
 5& 11 & 10 &8 \\
  &&&\\  \hline
 &&&\\
 9& 7 & 6 &12\\
  &&&\\  \hline
 &&&\\
4  &14&15& 1 \\
   &&&\\ \hline
   \end{tabular}
\end{table}
   \begin{table}
    \centering
   \caption{This matrix is a magic square for all real numbers
$a, b, c,$ and $d$ satisfying  (\ref{1}) and any real numbers
$k_1$ and $k_2$ }
    \label{tab5}
 \vspace{5mm}
\begin{tabular}{|c|c|c|c|}
 \hline
 \hspace{15mm} & \hspace{15mm} & \hspace{15mm} & \hspace{15mm} \\
$d+k_1+k_2$    &   $a$         &     $b$       & $c+k_1+k_2$ \\
  &&&\\    \hline
 &&&\\
$c+k_1$        & $b+k_2$       &   $a+k_2$     & $d+k_1$ \\
  &&&\\  \hline
 &&&\\
$c+k_2$        & $b+k_1$       &   $a+k_1$     & $d+k_2$ \\
  &&&\\  \hline
 &&&\\
$d$            &$a+k_1+k_2$    &   $b+k_1+k_2$ & $c$ \\
   &&&\\ \hline
   \end{tabular}
\end{table}

  \begin{table}[t]
    \centering
   \caption{The magic square with irrational entries obtained by using the data  (\ref{28}) }
    \label{tab6}
 \vspace{5mm}   {\small
\begin{tabular}{|c|c|c|c|}
 \hline
 \hspace{15mm} & \hspace{15mm} & \hspace{15mm} & \hspace{15mm} \\
$5+\sqrt{2}+\sqrt{3}$    &   $0$         &     $6$       & $1+\sqrt{2}+\sqrt{3}$ \\
  &&&\\    \hline
 &&&\\
$1+\sqrt{2}$        & $6+\sqrt{3}$       &   $\sqrt{3}$     & $5+\sqrt{2}$ \\
  &&&\\  \hline
 &&&\\
$1+\sqrt{3}$        & $6+\sqrt{2}$       &   $\sqrt{2}$     & $5+\sqrt{3}$ \\
  &&&\\  \hline
 &&&\\
$5$            &$\sqrt{2}+\sqrt{3}$    &   $6+\sqrt{2}+\sqrt{3}$ & $1$ \\
   &&&\\ \hline
   \end{tabular}   }
\end{table}

Let us now establish relations of the introduced method to the
theory of magic squares (and cubes) developed in
\cite{Adler1,Adler2}. A $T \times T$ matrix is called a {\it
magic square of order $T$} if   the sum of the entries in any
row, column, or diagonal is the same number. Usually, the entries
of the magic square are taken from the numbers $1, 2, \ldots ,
T^2$ and each of the  numbers is used exactly once as an entry.

There exist different ways to construct magic squares (see
\cite{Andrews}). Let us consider the one proposed in
\cite{Adler1,Adler2} and based on the Morse-Hedlund sequence. In
this paper we consider only the core case $T=4$ being a basis for
construction of the squares (and cubes) of the order $T=2^n$ (see
\cite{Adler1,Adler2} for details).

The algorithm to make the $4\times 4$ magic square is the
following. Count through the boxes in the square, and whenever
the number you count corresponds to `$0$' in the Morse-Hedlund
sequence (see Table~\ref{tab2}), write the number down in the
box. The rest of the boxes is filled in by counting numbers from
$16$ to $1$ and gives the magic square presented in
Table~\ref{tab4}.

Since the Prouhet identities were obtained as a particular case of
our method with the data from (\ref{27}) used in construction of
the magic square in Table~\ref{tab4}, let us execute the inverse
operation and substitute in Table~\ref{tab4} the numbers $1, 2,
\ldots , 16$ by the corresponding explicit formulae for the sets
$\{x_i: 1 \le i \le 8\}$ and $\{y_i: 1 \le i \le 8\}$. The result
is shown in Table~\ref{tab5} were, of course, the numbers $a, b,
c,$ and $d$ satisfy  (\ref{1}).

It can be shown that, due to (\ref{1}), the sum of  the entries
in any row, column, or diagonal is equal to $2(a+b+k_1+k_2)$.
Thus, the obtained matrix is a magic square for all real numbers
$a, b, c,$ and $d$ satisfying  (\ref{1}) and any real numbers
$k_1$ and $k_2$. Of course, some combination of these parameters
can give repeating numbers. Note that the data from (\ref{27})
give us the the magic square from Table~\ref{tab4} as a particular
case of Table~\ref{tab5}.

We conclude the paper by an example with irrational numbers.
 \beq
     a = 0,\hspace{3mm} b = 6,\hspace{3mm} c = 1,\hspace{3mm} d = 5, \hspace{1cm}
     k_1  = \sqrt{2}, \hspace{3mm}k_2  = \sqrt{3}.    \label{28}
\eeq
 For these data the corresponding systems (\ref{4}) and
(\ref{42}) become
\[
\hspace{-2cm}
 {\tiny   \begin{array}{rcl}
1^1  + 5^1  + \sqrt{2}^1  + (6+\sqrt{2})^1 + \sqrt{3}^1  + (6+\sqrt{3})^1  + (1+\sqrt{2}+\sqrt{3})^1  + (5+\sqrt{2}+\sqrt{3})^1  & = &  0^1  + 6^1  + (1+\sqrt{2})^1  + (5+\sqrt{2})^1  + (1+\sqrt{3})^1  + (5+\sqrt{3})^1  + (\sqrt{2}+\sqrt{3})^1  + (6+\sqrt{2}+\sqrt{3})^1,\\
1^2  + 5^2  + \sqrt{2}^2  + (6+\sqrt{2})^2 + \sqrt{3}^2  + (6+\sqrt{3})^2+ (1+\sqrt{2}+\sqrt{3})^2  + (5+\sqrt{2}+\sqrt{3})^2 & = &  0^2  + 6^2  + (1+\sqrt{2})^2  + (5+\sqrt{2})^2  + (1+\sqrt{3})^2  + (5+\sqrt{3})^2+ (\sqrt{2}+\sqrt{3})^2  + (6+\sqrt{2}+\sqrt{3})^2,\\
1^3  + 5^3  + \sqrt{2}^3  + (6+\sqrt{2})^3 + \sqrt{3}^3  +
(6+\sqrt{3})^3+ (1+\sqrt{2}+\sqrt{3})^3 + (5+\sqrt{2}+\sqrt{3})^3
& = &  0^3 + 6^3 + (1+\sqrt{2})^3  + (5+\sqrt{2})^3 +
(1+\sqrt{3})^3+ (5+\sqrt{3})^3+ (\sqrt{2}+\sqrt{3})^3  +
(6+\sqrt{2}+\sqrt{3})^3
\end{array}       }
\]
and
\[
\hspace{-2cm}
 {\tiny   \begin{array}{rcl}
1^1  + 5^1    & = &  0^1  + 6^1,\\
1^2  + 5^2  + \sqrt{2}^2  + (6+\sqrt{2})^2 & = &  0^2  + 6^2  + (1+\sqrt{2})^2  + (5+\sqrt{2})^2,\\
1^3  + 5^3  + \sqrt{2}^3  + (6+\sqrt{2})^3 + \sqrt{3}^3  +
(6+\sqrt{3})^3+ (1+\sqrt{2}+\sqrt{3})^3 + (5+\sqrt{2}+\sqrt{3})^3
& = &  0^3 + 6^3 + (1+\sqrt{2})^3  + (5+\sqrt{2})^3 +
(1+\sqrt{3})^3+ (5+\sqrt{3})^3+ (\sqrt{2}+\sqrt{3})^3  +
(6+\sqrt{2}+\sqrt{3})^3.
\end{array}       }
\]
Finally, the  magic square obtained from the data (\ref{28}) by
using the general method from Table~\ref{tab5} is shown in
Table~\ref{tab6}.

%{\small
%\[
%(a+k_1)^1+(b+k_1)^1=(c+k_1)^1+(d+k_1)^1
%\]
% \beq
% (a+k_1)^2+(b+k_1)^2+(c+k_1+k_2)^2+(d+k_1+k_2)^2=(c+k_1)^2+(d+k_1)^2+(a+k_1+k_2)^2+(b+k_1+k_2)^2
% \label{2}
% \eeq
%\[
% (a+k_1)^3+(b+k_1)^3+(c+k_1+k_2)^3+(d+k_1+k_2)^3
% +(c+k_1+k_3)^3+(d+k_1+k_3)^3+(a+k_1+k_2+k_3)^3+(b+k_1+k_2+k_3)^3
% =(c+k_1)^3+(d+k_1)^3+(a+k_1+k_2)^3+(b+k_1+k_2)^3
% +(a+k_1+k_3)^3+(b+k_1+k_3)^3+(c+k_1+k_2+k_3)^3+(d+k_1+k_2+k_3)^3
%\]}

\section{Appendix}

This Appendix contains an alternative description of the sets
$\{x_i: 1 \le i \le k\}$ and   $\{y_i: 1 \le i \le k\}$ proposed
by the unknown referee.

Given $n$, we construct $x_i, 1 \le i \le 2^{n+1},$ and   $y_i, 1
\le i \le 2^{n+1},$ so that for $1 \le m \le n+1$,
 \beq
\sum^{2^{n+1}}_{i=1}  x^m_i = \sum^{2^{n+1}}_{i=1}  y^m_i
             \label{a1}
 \eeq

Choose $a, b, c,$ and $d$ so that $a+b=c+d$. Choose any $n$ values
$k_1,  \ldots, k_{n}$. The $x_i$ are obtained by adding $a$ to an
even subset of the $k_j$ (including the empty set) or $b$ to an
even subset of the $k_j$ or $c$ to an odd subset of the $k_j$ or
$d$ to an odd subset of the $k_j$. The $y_i$ are obtained by
adding either $a$ or $b$ to an odd subset of the $k_j$ and $c$ or
$d$ to an even subset of the $k_j$.

Equation (\ref{a1}) is then equivalent to the statement that the
function of $a$ and $b$ given by
 \beq
f(a,b)=\sum_{S \subseteq \{ 1,2, \ldots, n \}}  (-1)^{|S|} \left[
\left(a +  \sum_{i \in S} k_i \right)^m +  \left(b +  \sum_{i \in
S} k_i \right)^m
       \right],      \label{a2}
 \eeq
is a function of $a+b$ when $1 \le m \le n+1$.

A simple expansion shows that
\begin{eqnarray*}
f(a,b) & = &  (a^m+b^m)\sum_{S} (-1)^{|S|} + \sum_{S \neq
\emptyset } (-1)^{|S|} \sum^{m-1}_{j=0} (a^j+b^j) \left(
\begin{array}{c} m
\\ j \end{array} \right) \left( \sum_{i \in S} k_i \right)^{m-j}\\
& = & \sum^{m-1}_{j=0} (a^j+b^j) \left(
\begin{array}{c} m
\\ j \end{array} \right) \sum_{S \neq
\emptyset } (-1)^{|S|} \left( \sum_{i \in S} k_i \right)^{m-j}\\
& = & \sum^{m-1}_{j=0} (a^j+b^j) \left(
\begin{array}{c} m
\\ j \end{array} \right) \sum_{\tiny{S = \{ s_1, \ldots, s_{|S|} \} \neq
\emptyset}}  (-1)^{|S|} \sum_{\tiny{i_1+ \ldots +i_{|S|}=m-j}}
\left( \small{\begin{array}{c} m-j
\\ i_1, \ldots, i_{|S|} \end{array}} \right) \prod^{|S|}_{j=1} k^{i_j}_{s_j}\\
& = & \sum^{m-1}_{j=0} (a^j+b^j) \left(
\begin{array}{c} m
\\ j \end{array} \right)  \sum_{i_1+ \ldots +i_n=m-j}  \left(
\begin{array}{c} m-j
\\ i_1, \ldots, i_n  \end{array} \right) \prod^{n}_{j=1} x^{i_j}_{j}
\sum_{S \ni  i_j > 0 \Rightarrow j \in S}  (-1)^{|S|}\\
& = & \sum^{m-1}_{j=0} (a^j+b^j) \left(
\begin{array}{c} m
\\ j \end{array} \right)  \sum_{\tiny{\begin{array}{c} i_1+ \ldots +i_n=m-j \\ \forall j,\,\,\, i_j \ge 1 \end{array}}}
\left( \begin{array}{c} m-j
\\ i_1, \ldots, i_n  \end{array} \right) \prod^{n}_{j=1} k^{i_j}_{j}
  (-1)^{n}\\
  & = & \left\{ \begin{array} {ll}
(-1)^{n}(n+1)!\,\,\, x_1 \cdots x_n (a+b+\sum^{n}_{i=1} k_{i}), &
\mbox{if }
m=n+1, \\
(-1)^{n}(n)!\,\,\, k_1 \cdots k_n, & \mbox{if } m=n, \\
0, & \mbox{if } m<n.
  \end{array} \right.
\end{eqnarray*}

 \vspace{5mm}
\begin{Acknowledgement}
The author thanks V.V. Korotkich for very useful and pleasant
discussions and the unknown referee for the alternative
description of the sets $\{x_i: 1 \le i \le k\}$ and   $\{y_i: 1
\le i \le k\}$ given in the Appendix.
\end{Acknowledgement}

%\end{article}

\begin{thebibliography}{}

\bibitem{Adler1}A. Adler (1992)  Magic cubes and 3-Adic Zeta function, The Mathematical
 Intelligencer,  14(3), 14--23.

\bibitem{Adler2}A. Adler  and S.-Y.R. Li (1977)  Magic cubes and Prouhet sequences,
American Mathematical Monthly,  84, 618--627.

\bibitem{Andrews}W.S. Andrews (1960) Magic squares and cubes,
Dover, New York.

\bibitem{[1]}   G.I. Arkhipov, A.A. Karatsuba and V.N.
Chubarikov  (1987) Theory  of  multiple
     trigonometric sums, Nauka, Moskow.

\bibitem{[15]}   M. Boldrin and M.  Woodford  (1990) Equlibrium  models  displaying
endogenous fluctuations and chaos: A survey, Journal  of  Monetary
     Economics, 25,
     189--220.

\bibitem{[10]}   G.J. Chaitin   (1966) On the length
of  programs  for  computing  finite binary
     sequences, J. Assoc. Comp. Math., 13, 547--570.


\bibitem{[13]}   M.J. Feigenbaum (1980) Universal behavior in  nonlinear  systems,
Los Alamos Science 1(1), 4--27.


\bibitem{[5]}    D. Hilbert (1909)  Beweis  f\"{u}r  die  Darstellbarkeit  der  ganzen
Zahlen durch eine feste Anzahl n-ter Potenzen, Math. Ann., 67,
      281--300.

\bibitem{Horst_and_Pardalos_(1995)}R. Horst and P.M. Pardalos (1995)   Handbook of Global Optimization, Kluwer Academic Publishers, Dordrecht.


\bibitem{[6]}    E. Kamke   (1921) Verallgeme inerungen des Waring - Hilbertschen  Satzes,
     Math. Ann., 3, 3--38.

\bibitem{[16]} A.A. Karatsuba  (1987) Hilbert-Kamke's problem in the analitical number
theory, Math. Zametki, 41(2), 272--284.



\bibitem{[11]}   A.N. Kolmogorov  (1965) Three approaches  to  the  quantitative
definition of
     information, Problems of Inform. Transmission, 1, 3--7.

\bibitem{[2]}    N.M. Korobov (1989) Trigonometric sums and their applications, Nauka,
     Moskow.

\bibitem{[8]}    V.V. Korotkich (1993) An  approach  to  nonlinear  systems  cycles
structural properties description, Institute for Automated
Systems, Moscow.

\bibitem{[9]}    V.V. Korotkich  (1999) A Mathematical Structure for Emergent
Computation, Kluwer Academic Publishers, Dordrecht.

\bibitem{[14]}   H.J. Korsch and H.-J. Jodl (1994) Chaos, Springer-Verlag,  Berlin.


\bibitem{Lehmer} D.H. Lehmer  (1947) The Tarry-Escott problem,
Scripta Math., 14, 37--41.

\bibitem{[3]}   L. Mordell   (1932) On a sum analogous to a Gauss's  sum,  Quart.  J.
Math., 3,  161--167.

\bibitem{Morse} M. Morse and G.A. Hedlund (1944) Unending chess,
symbolic dynamics and a problem in semigroups, Duke Math. J., 11,
1--7.

\bibitem{Prouhet}  E. Prouhet   (1851) Memoire sur quelques relations
 entre les puissances des nombres, C.R. Acad. Sci., Paris, 33, 225.

\bibitem{[12]}   R.J. Solomonoff  (1964) A formal theory of inductive inference,
Information and
     Control 7(1), 1--22.

\bibitem{Strongin_and_Sergeyev}R. G. Strongin and Ya. D.
Sergeyev (2000)    Global Optimization with Non-Convex
Constraints: Sequential and Parallel Algorithms, Kluwer Academic
Publishers, Dordrecht.


\bibitem{[4]}    R.C. Vaughan (1979) A survey of some important  problems  in  additive
     number theory, Societe Math. de France, Asterisque 61, 213--222.

\bibitem{[7]}    E. Waring (1770) Meditationes algebraicae, Cambridge.







\end{thebibliography}
\end{document}